\documentclass[11pt]{article}
	
\usepackage[margin=1in,marginparwidth=0.8in, marginparsep=0.1in]{geometry}
\usepackage{xy}
\usepackage{authblk}

\pdfoutput=1

\usepackage{soul}

\usepackage{amssymb,latexsym,amsmath,graphicx, bbm,tikz-cd}
\usepackage{stmaryrd} 
\usepackage{amscd}
\usepackage{cite}
\usepackage{todonotes}
\usepackage{color}
\usepackage{amsfonts} 

\usepackage[bookmarks=true, bookmarksopen=true,%
bookmarksdepth=3,bookmarksopenlevel=2,%
colorlinks=true,%
linkcolor=blue,%
citecolor=blue,%
filecolor=blue,%
menucolor=blue,%
urlcolor=blue]{hyperref}

\usepackage{times}
\usepackage{mathrsfs}
\usepackage{mathtools}
\usepackage{parskip}

\newtheorem{lemma}{Lemma}
\newtheorem{proposition}[lemma]{Proposition}
\newtheorem{corollary}[lemma]{Corollary}
\newtheorem{theorem}[lemma]{Theorem}

\newtheorem{mconstruction}[lemma]{Main Construction}

\newtheorem{definition}[lemma]{Definition}

\newtheorem{remark}[lemma]{Remark}
\newtheorem{example}[lemma]{Example}

 \newenvironment{proof}{{\bf
     Proof.}}{$\blacksquare$ \vspace{2mm}}

\newcommand{\A}{\mathbb{A}}

\newcommand{\G}{\mathbb{G}}

\newcommand{\N}{\mathbb{N}}

\renewcommand{\P}{\mathbb{P}}
\newcommand{\Q}{\mathbb{Q}}
\newcommand{\R}{\mathbb{R}}

\newcommand{\W}{\mathbb{W}}

\newcommand{\Z}{\mathbb{Z}}

\newcommand{\cD}{\mathcal{D}}

\newcommand{\cL}{\mathcal{L}}
\newcommand{\cM}{\mathcal{M}}
\newcommand{\cN}{\mathcal{N}}
\newcommand{\cO}{\mathcal{O}}

\newcommand{\Qb}{\mathbb{Q}}
\newcommand{\Zb}{\mathbb{Z}}
\newcommand{\Cb}{\mathbb{C}}
\newcommand{\Rb}{\mathbb{R}}
\newcommand{\Fb}{\mathbb{F}}
\newcommand{\Oo}{\mathcal{O}}
\newcommand{\Aa}{\mathcal{A}}

\newcommand{\Nn}{\mathcal{N}}

\newcommand{\m}{\mathfrak{m}}
\newcommand{\supp}{\mathsf{supp}}
\newcommand{\vol}{\mathsf{vol}}
\newcommand{\Spec}{\mathsf{Spec}\;}
\newcommand{\Arb}{\mathsf{Arb}}
\newcommand{\Xc}{\mathcal{X}}
\newcommand{\Uc}{\mathcal{U}}
\newcommand{\Spf}{\mathsf{Spf}\;}
\renewcommand{\sp}{sp}
\newcommand{\Dol}{\mathfrak{D}}
\newcommand{\pDol}{\mathfrak{D}}
\newcommand{\base}{B}
\newcommand{\pbase}{\Aa}
\newcommand{\pversal}{\mathfrak{V}}

\newcommand{\Kir}{\Psi}

\newcommand{\Ab}{\mathbb{A}}
\newcommand{\trop}{\mathsf{trop}}
\newcommand{\ctc}{\tau}
\newcommand{\periodmap}{\sigma}

\renewcommand{\top}{\mathsf{top}}

\newcommand{\redC}{\overline{C}}
\newcommand{\boing}{\mbox{$ \bigcirc \!\! \bullet$}}
\newcommand{\Bb}{\mathcal{{B}}}
\newcommand{\pcyclebase}{\Bb}
\newcommand{\hmap}{\pi_{\operatorname{res}}}
\newcommand{\hmapbasic}{\pi}
\newcommand{\Lpower}{N}
\newcommand{\rindex}{n}
\newcommand{\Wpower}{E}

\newcommand{\bmap}{\pi}

\newcommand{\stabset}{\operatorname{St}}
\newcommand{\gitbundle}{\cL_\eta} 
\newcommand{\gitbundlet}{\cL_{\tilde{\eta}}}
\newcommand{\edges}{E}
\newcommand{\edgesgam}{E(\Gamma)}
\newcommand{\bigtorus}{\mathbf{T}}
\newcommand{\diagtorus}{\G_m}
\newcommand{\AZ}{\A}
\newcommand{\fibera}{\hmap^{-1}(a)}
\newcommand{\unif}{\mathbf{u}}
\newcommand{\Kk}{R}

\begin{document}
\title{Hypertoric Hitchin systems and Kirchhoff polynomials}
\author[1]{Michael Groechenig\thanks{M.G. was supported by an NSERC discovery grant.}}
\author[2]{Michael McBreen\thanks{M.M. was supported by NSERC and the starting grant of Professor Artan Sheshmani at Aarhus University.}}
\affil[1]{Department of Mathematics, University of Toronto}
\affil[2]{ Department of Mathematics, Aarhus University and Harvard CMSA}
\maketitle
\abstract{We define a formal algebraic analogue of hypertoric Hitchin systems, whose complex-analytic counterparts were defined by Hausel--Proudfoot. These are algebraic completely integrable systems associated to a graph $\Gamma$. We study the variation of the Tamagawa number of the resulting family of abelian varieties, and show that it is described by the Kirchhoff polynomial of the graph $\Gamma$. In particular, this allows us to compute their $p$-adic volumes. We conclude the article by remarking that these spaces admit a volume preserving tropicalisation.}
\tableofcontents

\section{Introduction}

A \emph{network} $(\Gamma,w)$ is a graph $\Gamma = (V,E)$ (we permit multiple edges between two vertices) together with a weight function $w\colon E \to \Rb_{> 0}$. In his 1847 paper \cite{Kirchhoff}, Kirchhoff associated to a network the quantity
$$\Kir(\Gamma,w)=\sum_{T \in \Arb(\Gamma)} \prod_{e \notin T} w_e,$$ 
where $\Arb(\Gamma)$ is the set of maximal spanning trees of $\Gamma$ (respectively maximal spanning forests, if $\Gamma$ is disconnected). Curiously, the same polynomials also arise when computing Feynman integrals and their motivic counterparts, as in Bloch--Esnault--Kreimer\cite{BEK}. According to the main result of the present article, Kirchhoff polynomials compute fibrewise $p$-adic volumes of an abelian fibration associated to $\Gamma$. The complex-analytic counterpart of this fibration was constructed by Hausel--Proudfoot in an unpublished note, and has been studied in some detail in \cite{McBW} and \cite{DMS}. We use techniques similar to Mumford's \cite{mumford} to produce a formal algebraic counterpart of their construction.

\begin{mconstruction}[Hypertoric Hitchin systems]
For a regular Noetherian ring $R$ we denote by $\base(\Gamma)$ the affine scheme 
$$\Spec R[[T_e|e \in E]].$$ The union of coordinate hyperplanes in $\Ab_R^{E}$ defines a closed subset $\Delta \subset \base(\Gamma)$ whose set of irreducible components is in bijection with $E$. To a stability parameter $\eta\colon \edges \to \Rb/\Zb$ we associate a proper morphism
$$\pi\colon \pversal(\Gamma,\eta) \to \base(\Gamma),$$
such that 
$$\pi^{-1}(\base(\Gamma)\setminus \Delta) \to \base(\Gamma)\setminus \Delta$$
is a relative abelian scheme. Furthermore, for a generic choice of $\eta$, the total space $\pversal(\Gamma,\eta)$ is formally smooth. Under the same assumptions, the base change 
$$\Dol(\Gamma,\eta)=\pversal(\Gamma,\eta) \times_{\base(\Gamma)} (\base(\Gamma) \times_{\Ab^E_R} H^1(\Gamma,\Ab^1_R)) \to \base(\Gamma) \times_{\Ab^E_R} H^1(\Gamma,\Ab^1_R)$$
is an algebraic completely algebraic system. We refer to this morphism as the \emph{hypertoric Hitchin system}.
\end{mconstruction}

If $R = \Oo_F$ is the valuation ring of a local field (see \ref{sub:p-adic}) there is a canonical Borel measure on the topological space $\pversal(\Gamma,\eta)(\Oo_F)$. This is essentially the measure discussed in Andr\'e Weil's \cite{weil2012adeles}, although we need to slightly generalise his set-up (see Subsection \ref{sub:canonical} for details). This measure yields a notion of relative volume for the Hitchin fibres. Our main result asserts that this volume is computed in terms of the Kirchhoff polynomial.

\begin{theorem}[Main theorem]\label{main}
Let $b \in \base(\Gamma)(\Oo_F) \cap \left( \base(\Gamma) \setminus \Delta\right)(F)$, that is, an $\Oo_F$-rational point of the base, such that the induced $F$-rational point lies in the complement of $\Delta$. Then, the volume of the abelian $F$-variety $\pi^{-1}(b)$ is independent of the (generic) stability parameter $\eta$, and equals
$$\vol(\pi^{-1}(b)(F)) = (1-q^{-1})^{h^1(\Gamma)}\cdot{} \Kir(\Gamma,w_b),$$
where $w_b\colon E \to \Gamma$ is the function which sends an edge $e \in \Gamma$ to the value of the valuation $\nu_e(b)$, corresponding to the irreducible component of $\Delta$ associated to $e$.
\end{theorem}

The proof of our main result can be found in Corollary \ref{cor:main} in the main body of this text. Subsequently we turn to an interpretation of this computation.
In Section \ref{tropical} we show that these spaces admit a natural \emph{tropicalisation}. To be precise, we construct a piecewise linear manifold with corners $\trop(\pversal)(\Gamma,\eta)$ and a continuous map 
$$\trop(\pi)\colon \trop(\pversal)(\Gamma,\eta) \to [0,1)^E,$$
referred to as the \emph{tropical Hitchin map}. These spaces are endowed with a Borel measure. The following assertion is Proposition \ref{prop:trop}.

\begin{proposition}[Tropicalisation preserves volumes]
We retain the assumptions of Theorem \ref{main}. One has
$$\vol\left(\trop(\pi^{-1}(b))\right) = \vol\left(\pi^{-1}(b)\right).$$
\end{proposition}

The morphism $\pi$ resembles in many ways a relative compactified Jacobian of a relative nodal curve, yet it is not isomorphic to such a family. The last section is devoted to the observation that the tropicalisation of $\pi$ can be described as a Jacobian of a family of tropical curves. This suggests a connection between the theory of hypertoric Hitchin systems and recent work on tropicalisations of Jacobians \cite{baker-rabinoff,abreu-pacini,4authors}.

\begin{proposition}
Under the same hypotheses as in Theorem \ref{main}, the tropicalisation $\trop(\pi^{-1}(b))$ is isomorphic to $J(\Gamma,w_b)$ as a flat torus.
\end{proposition}

\bigskip
\noindent\textit{Acknowledgements.} The concept of hypertoric non-abelian Hodge theory, and thus in particular the hypertoric Hitchin map, is due to Tamas Hausel and Nick Proudfoot. Unfortunately their manuscript remains unpublished to this date, but without doubt their ideas permeate this article. We also gratefully acknowledge the influence of joint work of M. McB. with Vivek Shende and Zsuzsanna Dancso and of M. G. with Dimitri Wyss and Paul Ziegler. We thank B\'alint Vir\'ag for enlightening correspondence on Kirchhoff polynomials. It is a pleasure to thank everyone named above for their input and interesting conversations.

\section{An algebraic construction of hypertoric Hitchin systems}

In this section, we construct an algebraic counterpart of Dolbeault hypertoric spaces, together with their integrable system, which we will refer to as a hypertoric Hitchin system. These were originally defined by Hausel and Proudfoot in the complex analytic setting (unpublished), and studied from various perspectives in \cite{McBW} and \cite{DMS}. 
All schemes in this section are taken over $\Zb$. 

\subsection{Toric geometry of the universal cover of the Tate curve}
We start with infinitely many copies $\A^2_\rindex$ of the affine plane $\AZ^2$, indexed by $\rindex \in \Z$, with coordinates $x_\rindex, y_\rindex$. 
\begin{definition}\label{defi:W}
	Let $\W$ be the ind-scheme obtained by identifying $\AZ^2_\rindex \setminus \{x_\rindex = 0\}$ with $\AZ^2_{\rindex+1} \setminus \{y_{\rindex+1}=0\}$ via $y_{\rindex+1} = x_{\rindex}^{-1}$ and $x_\rindex y_\rindex = x_{\rindex+1}y_{\rindex+1}$.
	For a scheme $S$ we denote by $\W_S$ the base change $\W \times_{\Spec \Zb} S$.
\end{definition}
The map $x_\rindex y_\rindex \colon \A_n^2 \to \AZ^1$ extends by construction to a map $\hmapbasic : \W \to \AZ^1$. If we set $x = x_0$, $y=y_0$, we have the equality 
$$x_\rindex = \bmap^{\rindex}x \text{ and } y_\rindex = \bmap^{1-\rindex}x^{-1}$$ on $\W \setminus \bmap^{-1}(0)$. 

Let $R = \Z[\bmap, \bmap^{-1}]$; the basechange of $\W$ to $\Spec R$ is thus isomorphic to $\G_{m, R} = \Spec R[x,x^{-1}]$. On the other hand, the fibre $\W_0 := \bmap^{-1}(0)$ is an infinite chain of rational curves, glued end to end. Each chart $\A^2_{\rindex}$ intersects two such curves along its coordinate axes (the $x_{\rindex}$ and $y_\rindex$-axes). We denote the closure of the $y_{\rindex}$-axis by $\P^1_{\rindex}$. Note that this is also the closure of the $x_{\rindex+1}$-axis. We write $w_\rindex = \{ x_\rindex = y_\rindex = 0 \}$, and $(\G_m)_\rindex = \P^1_\rindex \setminus w_\rindex \cup w_{\rindex+1}$.

Let $\bigtorus = \G_{m} \times \G_{m}$, acting in the usual way on $\AZ^2$. This extends to an action on $\W$ preserving the charts $\A^2_\rindex$. The orbits of $\bigtorus$ are of three types:
\begin{align*}
 \W \setminus \W_0 & \cong \bigtorus  \\
 (\G_m)_\rindex & \ \ \  \rindex \in \Z \\
 w_\rindex & \ \ \ \rindex \in \Z.
\end{align*}
 $\W$ also carries an action of $\Z$ translating the copies of $\A^2_\rindex$, and commuting with $\bigtorus$. Its generator $1 \in \Z$ acts by $s_1^*x_{\rindex+1} = x_\rindex, s_1^* y_{\rindex+1} = y_\rindex$.

We define a line bundle $\cL$ on $\W$ by glueing trivial bundles on $\A^2_\rindex$ with transition function $x_\rindex$ on the overlaps $\A^2_\rindex \cap \A^2_{\rindex+1}$.
	\begin{remark}\label{rmk:twisty}The line bundle $\cL$ is naturally equivariant with respect to the $\bigtorus$ and $\Z$ actions. It is not, however, jointly equivariant. Instead, if $s_1^* \cL$ denotes the shift of $\cL$ by $1 \in \Z$, we have an equality of $\G_m$-equivariant bundles $s_1^*\cL = x\cL$ where we think of $x$ as a character of $\bigtorus$.
	\end{remark}
	
 In our given trivialisations, the section $f(x_\rindex, y_\rindex) \in H^0(\A^2_\rindex, \cL)$ is identified with the (rational) section $y_{\rindex+1} f(y_{\rindex+1}^{-1}, x_{\rindex+1}y_{\rindex+1}^2) \in H^0_{rat}(\A^2_{\rindex+1}, \cL)$.	
	
	It is sometimes more convenient to work with fixed coordinates $x, \pi$. We can restrict the given trivialisation of $\cL$ on $\A^2_\rindex$ to $\W \setminus \W_0 = \Spec \Z[x^{\pm 1}, \pi^{\pm 1}]$. We thus obtain an embedding $H^0(\W, \cL) \to \Z[x^{\pm 1}, \pi^{\pm 1}]$. The section $f(x, \pi) \in H^0(\A^2_{\rindex-1}, \cL)$ is identified with the rational section $\pi^{1-\rindex} x^{-1} f(x,\pi) \in H_{rat}^0(\A^2_{\rindex}, \cL)$.
\begin{proposition} \label{proposition:sections}
	The image of $H^0(\W, \cL) \to \Z[x^{\pm 1}, \pi^{\pm 1}]$ is the span of monomials $\bmap^kx^m$ such that $k \geq \frac{m^2 + m}{2}+1$. The section defined by a monomial is nonvanishing on $\W \setminus \W_0 \cup \G_{m} \cup \G_{m+1} \cup w_{m+1}$ if the inequality is an equality. If the inequality is strict, the section is nonvanishing only on $\W \setminus \W_0$.
\end{proposition}
	
\begin{proof}
A monomial $\pi^kx^m$ extends to a regular section of $\cL$ along $\A^2_\rindex$ if the expression $$\left( \prod_{i=1}^{\rindex} \pi^{1-i} x^{-1} \right) \bmap^k x^m = \left( \bmap^{-\frac{\rindex^2-\rindex}{2}} x^{-\rindex} \right) \bmap^k x^m = \bmap^{k - \frac{\rindex^2-\rindex}{2}}x^{m - \rindex}$$ is a polynomial in the chart coordinates $x_\rindex, y_\rindex$. 
Using the equality $$\bmap^Kx^M = \bmap^{K-\rindex M}x_\rindex^M = \bmap^{K + (1 - \rindex)M} y_n^{-M},$$ we find that this holds when
	\begin{align*}
	k + \frac{\rindex^2 - \rindex}{2} - m \rindex & \geq 0, \\ 
 k + \frac{\rindex^2 - 3\rindex}{2} -m \rindex  + m& \geq 0. 
 	\end{align*}
	The LHS of the first inequality is minimised at $n = m$ and $n=m+1$, where it equals $k - \frac{m^2+m}{2}$. 
	 The LHS of the second inequality is minimised at $n = m + 2$ and $n = m+1$, where it equals $k - \frac{m^2 + m}{2} - 1$.  The proposition follows.
\end{proof}

\begin{corollary} \label{cor:affinelocus}
	Let $s \in H^0(\W, \cL^\Lpower)$ be any nonzero $\bigtorus$-eigensection. Then the nonvanishing locus $\W_s$ of $s$ is isomorphic to either $\G^2_m, \G_m \times \A^1$ or $\A^2$. 
\end{corollary}

\begin{corollary}\label{cor:ample_mumford}
	The line bundle $\cL$ is ample in the sense used by Mumford \cite{mumford}, i.e., the nonvanishing loci of the sections of $\cL^n$ for $n \geq 1$ form a basis of the topology of $\W$. 
\end{corollary}

Every monomial $\bmap^k x^m$ defines a character $(k,m) \in \Z^2$ of $\bigtorus$. Let $\Delta \subset \R^2$ be the the unbounded polytope in $\R^2$ given by the convex hull of $ \mu_\rindex := (\frac{\rindex^2+\rindex}{2} + 1, \rindex)$ for $\rindex \in \Z$, and let $\Delta_\Z = \Delta \cap \Z^2$. 

By Proposition \ref{proposition:sections}, the $\bigtorus$-representation $H^0(\W, \cL)$ decomposes as a direct sum of one-dimensional representations $$H^0(\W, \cL) = \bigoplus_{\mu \in \Delta_\Z} H^0(\W, \cL)_{\mu}.$$

Fix $\Lpower > 1$. Let $\mu^N_\rindex := (\Lpower \frac{\rindex^2+\rindex}{2}+N, \Lpower \rindex)$ for $\rindex \in \Z$ and let $\Delta^\Lpower$ be the unbounded polytope in $\R^2$ given by the convex hull of $\mu^\Lpower_\rindex : \rindex \in \Z$. The same argument gives a decomposition into one-dimensional representations 
$$H^0(\W, \cL^\Lpower) = \bigoplus_{\mu \in \Delta^\Lpower_\Z} H^0(\W, \cL^\Lpower)_{\mu}.$$

\begin{lemma} \label{lem:supportbijection}
Let $S$ be a $\bigtorus$-orbit of $\W$. The characters of eigensections nonvanishing along $S$ are the integer points of a face of $\Delta^\Lpower$. This establishes a bijection between orbits and faces:
\begin{align*}
\W \setminus \W_0 & \to \Delta^\Lpower_\Z \\
(\G_m)_\rindex & \to [\mu^\Lpower_\rindex, \mu^\Lpower_{\rindex+1}] \\
w_\rindex & \to \Lpower \rindex.  
\end{align*}   
\end{lemma}
\begin{proof}
Follows from Proposition \ref{proposition:sections}.
\end{proof}

If we set $\Lpower=1$, there is still a map from orbits to faces defined by the same formulae, but it is no longer injective.

We have an antidiagonal torus $i : \G_m \to \bigtorus$ which preserves $\bmap$ and scales $x$ by the tautological character. 
\begin{lemma}  \label{stratumcharacters}
	 The induced map $i^* : \Z^2 \to \Z$ of character lattices defines surjections
\begin{align*}
\Delta^\Lpower_\Z & \to \Z \\
[\mu^\Lpower_\rindex, \mu^\Lpower_{\rindex+1}] & \to [Nn, N(n+1)] \\
\mu^\Lpower_{\rindex} & \to \Lpower \rindex.
\end{align*}    
\end{lemma}

\begin{definition} \label{def:orbitcharacters}
We write $\mu(S)$ for the set of $\diagtorus$-characters of sections of $\cL^\Lpower$ nonvanishing on $S$:
 \[ \mu(S) =
 \begin{cases}
 \Z \\
 [\Lpower \rindex, \Lpower (\rindex+1) ], & \rindex \in \Z \\
 \Lpower \rindex & \rindex \in \Z
 \end{cases}\] 
 \end{definition}

\subsection{Preliminaries on graphs}

Let $\Gamma$ be a finite oriented graph, i.e. a set of vertices $V(\Gamma)$, edges $E(\Gamma)$ and 'head and tail' maps $h: E(\Gamma) \to V(\Gamma), t : E(\Gamma) \to V(\Gamma)$ which record the vertices incident to an edge. We allow multiple edges between vertices, and edges which begin and end at the same vertex. We emphasise that the orientation is of purely auxiliary nature and will not affect the geometric properties of the spaces constructed below.

We will view $\Gamma$ as a CW complex in the natural way, so that for any commutative group scheme $\A$ we have natural identifications
\[ C_1(\Gamma, \A) = C^1(\Gamma, \A) = \A^{E(\Gamma)} \]
and
\[ C_0(\Gamma, \A) = C^0(\Gamma, \A) = \A^{V(\Gamma)}. \]
We write $\redC_0(\Gamma, \A)$ for the subgroup of chains summing to zero, and $\redC^0(\Gamma, \A)$ for the quotient by the constant chains.

We have the boundary and coboundary maps
\[ d_\Gamma : C_1(\Gamma, \A) \to \redC_0(\Gamma, \A) \]
and
\[ d_\Gamma^* : \redC^0(\Gamma, \A) \to C^1(\Gamma, \A). \]
An edge $e \in E(\Gamma)$ defines a chain in $C_1(\Gamma,\Zb)$. We will write $[e]$ for the class in $H^1(\Gamma,\Zb)$ which has value $1$ on this chain and $0$ on everything else. If $\gamma \in H_1(\Gamma,\Zb)$, we write $\langle \gamma, e \rangle$ for the natural pairing of $\gamma$ and $[e]$.

\subsection{The GIT quotient}
The group $\redC^1(\Gamma, \G_m) = \G_m^{\edgesgam}$ acts coordinatewise on $\W^{E(\Gamma)}$, where each factor acts by the antidiagonal embedding $\G_m \to \bigtorus$. Via the coboundary map $d^*_{\Gamma} $, this induces an action of $\redC^0(\Gamma, \G_m)$. Likewise, the group $C_1(\Gamma, \Z) = \Z^{\edgesgam}$ acts coordinatewise. This defines an action of $H_1(\Gamma, \Z)$ via the natural embedding. 

We define the line bundle $\cL^{E(\Gamma)}$ on $\W^{E(\Gamma)}$ as the box product of the bundles $\cL$ on the individual factors. It carries separate, non-commuting actions of $C^1(\Gamma, \G_m)$ and $C_1(\Gamma, \Z)$.  Given $y \in C_1(\Gamma,\Z)$, we have
\begin{equation} \label{eq:nonequivariance} y^* \gitbundle = y^\Lpower \otimes \gitbundle \end{equation}
In particular, the image of the coboundary map $\redC^0(\Gamma, \G_m) \to C^1(\Gamma, \G_m)$ commutes with the action of $H_1(\Gamma,\Z)$ on $\gitbundle$. 

Given $\eta \in \redC_0(\Gamma, \Z)$, and $N \gg 1$, we write $\gitbundle$ for the $\redC^0(\Gamma, \G_m)$-equivariant bundle obtained by taking the $\Lpower$-th power of $\cL^{E(\Gamma)}$ and twisting by the character $\eta$.

\begin{definition}
Fix $\eta, N$ as above.
\begin{itemize}
	\item Let $\widetilde{\pversal}(\Gamma, \eta) := \W^{E(\Gamma)} \sslash_{\gitbundle} \redC^0(\Gamma, \G_m)$ be the quotient in the sense of GIT with linearisation $\gitbundle$. It carries commuting actions of $H^1(\Gamma, \G_m)$ and $H_1(\Gamma, \Z)$.
	
	\item Let $\cL^\Lpower_\eta$ be the descent of the linearisation bundle; it is a $C^1(\Gamma, \G_m)/\redC^0(\Gamma, \G_m) = H^1(\Gamma, \G_m)$-equivariant line bundle on the quotient.
	
	\item Let $\hmap \colon \widetilde{\pversal}(\Gamma, \eta) \to C_1(\Gamma, \A^1)$ be the map defined by descending $\pi^{\edgesgam} \colon \W^{\edgesgam} \to C_1(\Gamma, \A^1)$ to the quotient.
	\end{itemize}
\end{definition}

To study this GIT quotient, we first consider the $C^1(\Gamma, \bigtorus)$-orbits of $\W^{E(\Gamma)}$. Any such orbit is a product\[ S = \prod_{e \in E(\Gamma)} S(e)\] of $\bigtorus$-orbits on $\W$ indexed by $e \in \edgesgam$, where $S(e)$ can be any of
\[ S(e) =
\begin{cases}
\W \setminus \W_0 \\
(\G_m)_\rindex, & \rindex \in \Z \\
w_\rindex & \rindex \in \Z
\end{cases}\] 

The closure of these orbits are denoted by  
\[ \overline{S}(e) =
\begin{cases}
\W \\
\P^1_\rindex, & \rindex \in \Z \\
w_\rindex & \rindex \in \Z
\end{cases}\]

Let $R(S)$ be the representation of $C^1(\Gamma, \G_m)$ defined by the image of the restriction $H^0(\W^\Wpower, \cL^\Lpower) \to H^0(S, \cL^\Lpower)$. 
\begin{proposition}
	We have a decomposition into $C^1(\Gamma, \G_m)$-isotypic components 
	$$R(S) = \bigoplus_{\mu \in \mu(S)} R(S)_\mu$$ where
	\begin{equation}
	\mu(S) := \prod_{e \in E(\Gamma)} \mu(S(e)) \subset \Z^{E(\Gamma)} 
	\end{equation}
	as in Definition \ref{def:orbitcharacters} and $\Z^{E(\Gamma)} = C_1(\Gamma, \Z)$ is the character lattice of $C^1(\Gamma, \G_m)$. 
\end{proposition}

\begin{remark} Note that $R(S)_\mu$ will typically not be one-dimensional. The exception is the case where $S$ contains a dense orbit of $C^1(\Gamma, \G_m)$, or equivalently, $S \subset \W_0^{\Wpower}$. \end{remark}

 Let $R(S, \eta)$ be the representation of $\redC^0(\Gamma, \G_m)$ defined by the image of the restriction $H^0(\W^\Wpower, \gitbundle) \to H^0(S, \gitbundle)$. We wish to describe its invariant part. To do so it is useful to lift it to a representation of $C^1(\Gamma, \G_m)$. 
 
 Thus, lift $\eta$ to a character of $\tilde{\eta}$ of $C^1(\Gamma, \G_m)$ and let $R(S, \tilde{\eta})$ be the representation of $C^1(\Gamma, \G_m)$ defined by the image of the restriction $H^0(\W^\Wpower, \gitbundlet) \to H^0(S, \gitbundlet)$. 

\begin{lemma} \label{lem:components}
We have a decomposition into $C^1(\Gamma, \G_m)$-isotypic components $$R(S, \tilde{\eta})^{\redC^0(\Gamma, \G_m)} = \bigoplus_{\mu} R(S, \tilde{\eta})_\mu$$ 
where the component is nontrivial if and only if $$\mu \in d_{\Gamma}^{-1}(0) \cap (\mu(S) +\tilde{\eta}).$$
\end{lemma}
 Translation by $\tilde{\eta}$ identifies this intersection with $d_{\Gamma}^{-1}(-\eta) \cap \mu(S)$, which as expected does not depend on the choice of lift $\tilde{\eta}$. Recall that a geometric point of $\W^{E(\Gamma)}$ is semistable for  the action of $\redC^0(\Gamma, \G_m)$ linearised by $\gitbundle$ if there exists a $\redC^0(\Gamma, \G_m)$-invariant section $H^0(\W^{E(\Gamma)}, \gitbundle)$ nonvanishing at this point, with affine nonvanishing locus.

Since $\gitbundle$ is $C^1(\Gamma, \bigtorus)$-equivariant, it is not hard to see that if one point in $S$ is semistable, then so are all other points in $S$.
\begin{lemma}
The orbit $S$ is semistable if and only if $d_{\Gamma}^{-1}(-\eta) \cap \mu(S)$ is nonempty.
\end{lemma}
\begin{proof}
By Lemma \ref{lem:components}, this is equivalent to the existence of a nonvanishing invariant section, which we may choose to be an eigensection of the bigger torus $C^1(\Gamma, \bigtorus)$. Corollary \ref{cor:affinelocus} shows that such a section will have affine nonvanishing locus.
\end{proof}

We now turn to the question of stabilisers. The only fixed points of the antidiagonal $\G_m$-action on $\W^2$ are $w_\rindex, \rindex \in \Z$. Let $E(S,w) \subset E(\Gamma)$ be the subset of edges for which $S(e) = w_\rindex$ for some $\rindex$. The stabiliser of $S$ in $\redC^0(\Gamma, \G_m)$ is the kernel of $$d^*_{\Gamma \setminus E(S, w)} : \redC^0(\Gamma, \G_m) \to C^1(\Gamma \setminus E(S, w), \G_m).$$ In particular, $\redC^0(\Gamma, \G_m)$ acts freely on $S$ if and only if $\Gamma \setminus E(S, w)$ is connected.
\begin{definition}
 We say $\eta \in C_0(\Gamma, \Z)$ is generic if $\Gamma \setminus E(S, w)$ is connected for all $S$ for which  $d_{\Gamma}^{-1}(-\eta) \cap \mu(S)$ is nonempty.
\end{definition}
It follows that the $\redC^0(\Gamma, \G_m)$-action on the semistable locus is free for generic $\eta$. 
\begin{corollary}
	For generic $\eta$, the semistable and stable loci coincide, and $\widetilde{\pversal}(\Gamma, \eta)$ is smooth.
\end{corollary}

\subsection{The fibers of $\hmap$} \label{sec:fibersofpi}
Fix $\eta$, not necessarily generic. We will now describe the fibers of $\hmap$ over the locus in $C_1(\Gamma, \Ab^1)$ where none of the edge coordinates vanish, i.e. $C_1(\Gamma, \G_m)$. Let $\Kk_{\Gamma} = \Z[\pi_e, \pi_e^{-1}]_{e \in \edgesgam}$, so that $C_1(\Gamma, \G_m) = \Spec \Kk_{\Gamma}$. By definition, we have 
\begin{align*}
\widetilde{\pversal}^{alg}(\Gamma, \eta) \times_{C_1(\Gamma, \Ab^1)}C_1(\Gamma, \G_m)  = (\W \setminus \W_0)^{\edgesgam} \sslash_{\gitbundle} \redC^0(\Gamma, \G_m) \end{align*}
Since all $\redC^0(\Gamma, \G_m)$ orbits in $ (\W \setminus \W_0)^{\edgesgam}$ are stable and $(\W \setminus \W_0) \cong \G_{m, K}$, we have 
\[  (\W \setminus \W_0)^{\edgesgam} \sslash_{\gitbundle} \redC^0(\Gamma, \G_m) = C^1(\Gamma, \G_{m,\Kk_{\Gamma}}) / \redC^0(\Gamma, \G_{m,\Kk_{\Gamma}}) \]  
Collecting the above, we obtain the following description of the generic fibers.
\begin{lemma} \label{lem:genericfiber}
 \[ \widetilde{\pversal}^{alg}(\Gamma, \eta) \times_{C_1(\Gamma, \Ab^1)}C_1(\Gamma, \G_m) = H^1(\Gamma, \G_{m,\Kk_{\Gamma}}). \] 
 
 The action of $H_1(\Gamma, \Z)$ on this locus is given by the composition
\begin{equation} \label{eq:periodlattice} \sigma \colon H_1(\Gamma, \Z) \to C^1(\Gamma, \G_{m,\Kk_{\Gamma}}) \to H^1(\Gamma, \G_{m,\Kk_{\Gamma}}) \end{equation} where the first map is 
\[ \gamma \to (\pi^{\langle \gamma, e \rangle}_e )_{e \in \edgesgam} \]
and the second map is the natural projection.
\end{lemma}

	We now turn to the central fiber $\hmap^{-1}(0) = \W_0^{E(\Gamma)} \sslash_{\gitbundle} \redC^0(\Gamma, \G_m)$. The ind-scheme $\W_0^{E(\Gamma)}$ is the union of $C^1(\Gamma, \bigtorus)$-orbits $S$ with 
	\[ S(e) = \begin{cases}
	(\G_m)_\rindex \\
	w_\rindex
	\end{cases}. \] 
	Those $S$ for which $S(e) = (\G_m)_{\rindex_e}$ (equivalently, $\overline{S}(e) = \P^1_{\rindex_e}$) for all $e \in \Gamma$ are free orbits of $C^1(\Gamma, \G_m)$. We write $\vec{n} \in \Z^{E(\Gamma)}$ for the tuple of values of $\rindex_e$, and label the orbit $S_{\vec{\rindex}}$. The closures $\overline{S}_{\vec{\rindex}}$ are the irreducible components of $\W_0^{E(\Gamma)}$, isomorphic to $(\P^1)^{E(\Gamma)}.$ These components intersect along the closures of smaller orbits, isomorphic to $(\P_{\Zb}^1)^{k}$ for $k < |E(\Gamma)|.$ 	
	
	The GIT quotient of a union of proper subvarieties stable under the acting group is the union of the GIT quotients. We thus obtain the following picture of the fibre.
	
\begin{lemma} \label{lemma:specialfibreunion}
	 The special fibre $\hmap^{-1}(0)$ is the union of subvarieties $$\frak{X}_{\vec{\rindex}} := \left( \prod_{e \in E(\Gamma)} \P^1_{\rindex_e} \right) \sslash_{\gitbundle} \redC^0(\Gamma, \G_m)$$ over all $\vec{\rindex} \in \Z^n$ such that 
	\begin{equation} \label{eq:existencecondition}  \mu(S_{\vec{\rindex}})  \cap d_{\Gamma}^{-1}(-\eta) = \left( \prod_{e \in E(\Gamma)} [N \rindex_e, N(\rindex_e+1)] \right) \cap d_{\Gamma}^{-1}(-\eta) \neq \emptyset. \end{equation}
	The intersection $\frak{X}_{\vec{\rindex}} \cap \frak{X}_{\vec{\rindex}'}$ is nonempty if and only if 
	    	\begin{equation} \label{eq:existencecondition2} d_{\Gamma}^{-1}(-\eta) \cap \mu(S_{\vec{\rindex}}) \cap  \mu(S_{\vec{\rindex}'}) \neq \emptyset. \end{equation}
	in which case $\vec{\rindex}$ and $\vec{\rindex}'$ differ by by at most one in each coordinate, and we have
	   $$\frak{X}_{\vec{\rindex}} \cap \frak{X}_{\vec{\rindex}'} = \left( \prod_{\rindex_e = \rindex'_e+1} w_{\rindex_{e}} \times \prod_{\rindex_e = \rindex'_e-1} w_{\rindex'_{e}} \times \prod_{\rindex_e = \rindex'_e} \P^1_{\rindex_e} \right) \sslash_{\cL^\Lpower \otimes \eta} \redC^0(\Gamma, \G_m).$$ 
\end{lemma}
\begin{corollary} \label{cor:specialfibreconnectedandproper}
	The special fibre $\hmap^{-1}(0)$ is connected, and its irreducible components are proper toric varieties.
\end{corollary}
\begin{proof}
	Properness is clear from Lemma \ref{lemma:specialfibreunion}. Connectedness follows from the fact that $\mu(S_{\vec{n}})$ tile $\Z^{E(\Gamma)}$.
\end{proof}

Let $\stabset_\eta \subset \Z^{E(\Gamma)} = C_1(\Gamma, \Z)$ be the set of $\vec{\rindex}$ satisfying \eqref{eq:existencecondition}. Consider the action by translations of $H_1(\Gamma, \Z)$ on $C_1(\Gamma, \Z)$. 
\begin{lemma} \label{lem:finiteorbitset}
	For any $N>0$, there are finitely many $H_1(\Gamma,N\Z)$-orbits in $\stabset_\eta$.
\end{lemma}
\begin{proof}
By construction, $d_\Gamma^{-1}(\eta)_\R$ is a torsor over $H_1(\Gamma, \R)$. It follows that $d_\Gamma^{-1}(\eta)_\R / H_1(\Gamma, N\Z)$ is compact. The tiling of $C_1(\Gamma, \R)$ by the cubes $\mu(S_{\vec{n}})_\R$ defines a stratification whose strata are interiors of cube faces, with the interiors of the cubes as maximal strata. By compactness, $d_\Gamma^{-1}(\eta)_\R / H_1(\Gamma, N\Z)$ is covered by the image of finitely many of these strata.
\end{proof}
\begin{corollary} \label{cor:finitetypetiles}
Fix any integer $N > 0$. There is an open subset $U \subset \widetilde{\pversal}^{alg}(\Gamma, \eta)$, of finite type over $C_1(\Gamma, \A)$, such that the union of $y \cdot U$ for $y \in H_1(\Gamma, N\Z)$ covers $\widetilde{\pversal}^{alg}(\Gamma, \eta)$.
 \end{corollary}
 \begin{proof}
For any $\vec{\rindex} \in C_1(\Gamma, \Z)$, let $\square_{\vec{\rindex}} \subset C_1(\Gamma, \Z)$ be the set of lattice points in a box of sidelength $2N$ centered at $\vec{\rindex}$. Let $U_{\vec{\rindex}}$ be the union over $\vec{\rindex}' \in \square_{\vec{\rindex}}$ of $\prod_{e \in \edgesgam} \A^2_{\rindex'} \subset \W^{\edgesgam}$. By Lemma \ref{lem:finiteorbitset}, there is a finite set $F \subset \stabset_\eta$ whose $H_1(\Gamma,N\Z)$-orbits cover $\stabset_\eta$. Let $\hat{U}$ be the union of $U_{\vec{\rindex}}$ for $\vec{\rindex}$ in $F$. The image of $\hat{U}$ in $ \widetilde{\pversal}^{alg}(\Gamma, \eta)$ is the requisite subset.
 \end{proof}

\subsection{Mumford's approach to quotients}

In this subsection we explain how to ``quotient" the ind-scheme $\widetilde{\pversal}^{alg}(\Gamma, \eta)$ by the discrete group $H^1(\Gamma,\Zb)$ in an algebraic manner. The inverted commas serve as reminders that we first have to replace $\widetilde{\pversal}^{alg}(\Gamma, \eta)$ by the formal completion of its central fibre, then take the quotient, and then \emph{algebraise} the resulting formal scheme to ${\pversal}^{alg}(\Gamma, \eta)$. 

The ``quotient" procedure outlined above is due to Mumford (see \cite{mumford}). It is difficult to improve upon the exposition of the aforementioned article, and we therefore only give a rough overview of his ideas.

Let $A$ be an excellent integrally closed noetherian ring, $I \subset A$ an ideal with $\sqrt{I} = I$, such that $A$ is complete with respect to the $I$-adic topology. Let $K$ be the quotient field of $A$. Mumford \cite{mumford} gives a recipe for constructing semi-abelian schemes over $\Spec A$, whose generic fibre is an abelian variety and whose fibre over $\Spec A/I$ is $\G_{m,A/I}^r$. Following Mumford, we will denote the group scheme $\G_{m,A}^r$ by $\tilde{G}$.

The requisite input data is a lattice embedding $Y \to \tilde{G}(K)$ and a `polarisation' $\phi : Y \to X$ where $X$ is the character lattice of $\tilde{G}$. Given $x \in X$, let $\mathcal{X}^x : \tilde{G}(K) \to K$ be the associated function. The data is required to satisfy:
\begin{enumerate}
\item $\mathcal{X}^{\phi(y)}(z) = \mathcal{X}^{\phi(z)}(y)$ for all $y,z \in Y$.
\item $\mathcal{X}^{\phi(y)}(y) \in I$ for all $y \in Y$.
\end{enumerate}

\begin{definition}\label{defi:models}
A \emph{relatively complete model} of $(\tilde{G},Y,\phi)$ is given by 
\begin{enumerate}
\item[(a)] an integral $A$-scheme $\tilde{P}$, locally of finite type, endowed with actions by $\tilde{G}$ and $Y$;
\item[(b)] an open immersion $i\colon \tilde{G} \hookrightarrow \tilde{P}$ equivariant with respect to $\tilde{G}$ and $Y$;
\item[(c)] a line bundle $\tilde{L}$ on $\tilde{P}$, equivariant with respect to $\tilde{G}$ and $Y$.
\end{enumerate}
This data must satisfy the following conditions:
\begin{enumerate}
\item There is an open subset $U \subset \tilde{P}$, which is $\tilde{G}$-invariant and of finite type over $A$, such that $\tilde{P}$ is covered by the union of all $Y$-translates of $U$
$$\tilde{P} = \bigcup_{y \in Y} y\cdot{} U.$$

\item The line bundle $\tilde{L}$ is \emph{ample} in the sense that the open subsets 
$$V(s) = \{x \in \tilde{P}| s(x)\text{ invertible}\}$$
form a basis for the topology of $\tilde{P}$, where $s$ ranges through sections of $\tilde{L}^n$ for $n$ a positive integer.

\item Let $S_y\colon \tilde{L} \to \tilde{L}$ denote the action of $y \in Y$ on total space of $\tilde{L}$. Similarly, for every $A$-scheme $S$, and an $S$-valued point $a \in \tilde{G}(S)$, let $T_a\colon \tilde{L}\times_{\Spec A} S \to \tilde{L} \times_{\Spec A} S$ be the analogous map obtained from the $\tilde{G}$-equivariant structure. Then we have
$$T_a S_y = \mathcal{X}^{\phi(y)}(a)\cdot{}S_y T_a$$
on $\tilde{L} \times_{\Spec A} S$

\item Let $\tilde{P}_0$ denote the base change $\tilde{P} \times_{\Spec A} \Spec A/I$. Then $\tilde{P}_0$ is connected and every irreducible component of $\tilde{P}_0$ is proper over $\Spec A/I$.

\end{enumerate}
\end{definition}

\begin{remark}
The definition is not exactly the same as \cite[Definition 2.1]{mumford}. Instead of the fourth condition, Mumford requires a stronger assumption of valuation-theoretic nature. The latter implies the assumptions above on $\tilde{P}_0$ (see Proposition 3.3 and Theorem 3.8 in \emph{loc. cit.}). The proof of Mumford's Theorem 3.10 (see Theorem \ref{thm:mumford} below) only requires connectedness of $\tilde{P}_0$ and properness of irreducible components of $\tilde{P}_0$.
\end{remark}

Section 3 of \cite{mumford} is devoted to the proof of the following result (Theorem 3.10 in \emph{loc. cit.}), which constructs a ``quotient" of $\tilde{P}$ by $Y$. We denote by $\mathfrak{m} \subset A$ the maximal ideal, and the fibre product
$$\tilde{P} \times_{\Spec A} \Spec A/\mathfrak{m}^n$$
by $\tilde{P}_n$.

\begin{theorem}[Mumford]\label{thm:mumford}
Assuming that the assumptions of Definition \ref{defi:models} hold, there exists a family of projective $A/\mathfrak{m}^n$-schemes $P_n$, indexed by $n \in \mathbb{N}$, together with an ample line bundle $L_n$, such that there is an \'etale surjective morphism
$$\rho\colon\tilde{P}_n \to P_n,$$
satisfying the conditions
\begin{enumerate} 
\item $\rho^*L_n \simeq \tilde{L}|_{\tilde{P}_n}$,
\item $\rho(x) = \rho(y)$ if and only $x,y \in \tilde{P}_n$ lie in the same $Y$-orbit.
\end{enumerate}
\end{theorem}

By virtue of the formal existence theorem, Mumford therefore obtains the existence of a projective $A$-scheme $P$, such that $P_n = P \times_{\Spec A} \Spec A/\mathfrak{m}^n$ for all $n \in \mathbb{N}$.

Henceforth we fix a regular Noetherian ring $R$ and denote it spectrum by $S = \Spec R$. All schemes and morphisms will be $S$-schemes. This step is necessary since we will work with schemes obtained by completions of rings. This is a process which does not commute with base change, and hence we are forced to change to $S$-schemes before completions enter the picture.

We denote by $A$ the completion of the local ring of $C_1(\Gamma,\Ab^1_S)$ at $0$.  Note that there is a natural map $\Spec A \to C_1(\Gamma,\Ab^1_{S})$. 
\begin{proposition}
There exists a line bundle $\tilde{L}$ on the base change 
$$\tilde{P}=\tilde{\pversal}(\Gamma,\eta) \times_{C_1(\Gamma,\Ab^1_{S})} \Spec A,$$
such that $(\tilde{P},\tilde{L},\phi)$ is a relatively complete model for $\tilde{T}=H^1(\Gamma,\G_{m,S})$.
\end{proposition}

\begin{proof}
It remains to define $\tilde{L}$ and to verify the defining properties for relatively complete models. Since $\tilde{\pversal}(\Gamma,\eta)$ is defined as a GIT quotient, it is endowed with an ample line bundle $\tilde{L}$, induced by $\gitbundle$ on $\W^{\edgesgam}$. Let us check that the triple $(\tilde{P},\tilde{L},\phi)$ satisfies the conditions of Definition \ref{defi:models}. \begin{enumerate}
\item This follows from Corollary \ref{cor:finitetypetiles}. 
\item This follows from ampleness of $\cL$ established in Corollary \ref{cor:ample_mumford}.
\item This follows from the analogous property for $\cL$ that we recorded in Remark \ref{rmk:twisty}.
\item By virtue of Corollary \ref{cor:specialfibreconnectedandproper}, the irreducible components of the special fibre $\tilde{P}_0$ are proper toric schemes.
\end{enumerate}
This concludes the proof of the proposition.
\end{proof}

\begin{definition} \label{def:versalspace}
Let $\pversal(\Gamma,\eta)$ be the projective $A$-scheme obtained from the relatively complete model $(\tilde{P},\tilde{L},\phi)$. We write $\hmap \colon \pversal(\Gamma,\eta) \to \Spec A$ for the natural map.
\end{definition}
By construction, the base change
$$\pversal(\Gamma,\eta) \times_{C_1(\Gamma, \A^1_S)} C_1(\Gamma,\G_{m,S})$$
is an abelian scheme.

\begin{definition}
We define $\cD(\Gamma,\eta)=\pversal(\Gamma,\eta) \times_{C_1(\Gamma,\Ab^1_{S})} H_1(\Gamma,\Ab^1_S)$.
\end{definition}
This is an algebraic analogue of the hypertoric Dolbeault space. As in the complex-analytic category, one can show that for a generic stability parameter $\eta$, the scheme $\cD(\Gamma,\eta)$ is endowed with a (formal) symplectic structure, and the morphism $\hmap \colon \cD(\Gamma,\eta) \to \Spec A \times_{C_1(\Gamma,\Ab^1_S)} H_1(\Gamma,\Ab^1_S)$ is an algebraic integrable system (i.e., the generic fibre is Lagrangian).

\begin{remark}In \cite[Theorem 2.5]{mumford} Mumford shows that relatively complete models always exist after replacing $\phi$ by an integral power. His construction is different from the one outlined above. It relies on the choice of a finite subset $\Sigma \subset X$ containing $0$, stable under $x \to -x$ and containing a basis. Mumford calls such a set a {\em star}. One then defines $\frak{P}$ by taking the projectivisation of an explicit graded ring $R_\Sigma$. 
Presently we do not know how to relate the choice of a stability condition with Mumford's choice of a star. 
\end{remark}

\section{Volume computation}

\subsection{Preliminaries on $p$-adic integration}\label{sub:p-adic}

We denote by $F$ a non-archimedean local field. It contains a unique compact-open subring $\Oo_F \subset F$, which is a local ring. Furthermore, we will use the following notation ...
\begin{itemize}
\item[...] $\m_F \subset \Oo_F$ for the unique maximal ideal, 
\item[...] $k_F$ for the residue field,
\item[...] $q=|k_F|$,
\item[...] $p=\mathsf{char}(k_F)$.
\end{itemize}
The field of $p$-adic numbers $\Qb_p$ is a non-archimedean local field. One has
$$\Oo_{\Qb_p} = \Zb_p\text{  and  }k_{\Qb_p} = \Fb_p.$$
Any finite extension of $\Qb_p$ is also a non-archimedean local field. The field of formal Laurent series
$\Fb_q((t))$ is another example. Its ring of integers and residue field are 
$$\Oo_{\Fb_q((t))} = \Fb_q[[t]] \text{,  respectively  }k_{\Fb_q((t))}=\Fb_q.$$

Despite being totally disconnected, these fields share many features with real and complex numbers. For instance, there is a well-behaved theory of $F$-analytic manifolds:

\begin{definition}
\begin{enumerate}
\item[(a)] Let $U \subset F^m$ be an open subset of a finite-dimensional $F$-vector space. A function
$$f\colon U \to F^n$$
is said to be \emph{analytic} it is locally expressible by a converging power series.

\item[(b)] Let $M$ be a topological space. A set of pairs $\mathcal{U}=\{(U,\phi)\}$ consisting of open subsets $U \subset M$ and homeomorphisms
$$\phi\colon U \xrightarrow{\simeq} U',$$
where $U'$ is an open subset of a finite-dimensional $F$-vector space $F^n$, is said to be a \emph{maximal $F$-atlas} if the following conditions are met:
\begin{itemize}
\item For $(U_i,\phi_i)_{i=1,2} \in \mathcal{U}$ the associated change-of-coordinates map
$$\phi_2 \circ \phi_1^{-1}\colon \phi_1(U_1) \to \phi_2(U_2)$$
is analytic.
\item A pair $(V,\psi)$ consisting of an open subset $V \subset M$ and a homeomorphism $\psi\colon V \xrightarrow{\simeq} V'$ belongs to $\mathcal{U}$, if and only if for every $(U,\phi) \in \mathcal{U}$ the associated change-of-coordinates map
$$\psi \circ \phi^{-1}\colon \phi(U) \to \psi(V)$$
is analytic.
\end{itemize}

\item[(c)] A topological space $M$ endowed with a maximal $F$-atlas $\mathcal{U}$ is said to be an \emph{$F$-analytic manifold}, if is Hausdorff and second-countable.
\end{enumerate}
\end{definition}

Similar to the theories of real and complex manifolds one can define the notions of analytic functions, vector fields, $1$-forms, and more generally $i$-forms. We leave the details to the reader, and only recall the formal definition of an analytic top-degree form.

\begin{definition}
Let $M$ be an $F$-analytic manifold. A rule which assigns to every chart $(U,\phi) \in \mathcal{U}$ an analytic function
$$\omega_{(U,\phi)}\colon \phi(U) \to F$$
is said to be an \emph{analytic top-degree form}, if for every pair of charts $(U_i,\phi_i)_{i=1,2} \in \mathcal{U}$ we have an equality
$$\omega_{(U_2,\phi_2)} = \det\left( \partial(\phi_2\circ \phi_1^{-1}) \right) \cdot{} \omega_{(U_1,\phi_1)}$$
over the open subset $\phi_2(U_1\cap U_2) \subset \phi_2(U_2)$.
\end{definition}

By virtue of assumption, a finite-dimensional $F$-vector space $F^n$ is a locally-compact abelian group (with respect to addition). We denote by $\mu$ the Haar measure on $F^n$ satisfying $\mu(\Oo_F^n) = 1$.

Let $\omega$ be an analytic top-degree form on $M$, such that $\supp(\omega)$ is contained in a single chart $(U,\phi)$. The integral

\begin{equation}\label{eqn:1}\int_M |\omega| = \int_{\phi(U)} |\omega_{(U,\phi)}| d\mu\end{equation}
is independent of the chosen chart $(U,\phi) \supset \supp(\omega)$. Using partitions of unity, one can extend this definition to arbitrary analytic top-degree forms.

\begin{definition}
We denote by $\omega\mapsto \int_M |\omega|$ the unique function, additive in $\omega$, which agrees with \eqref{eqn:1} for top-degree forms supported within a single chart.
\end{definition}

We emphasise that the integral above is an integral of a real-valued function, since we integrate the absolute value $|\omega|$ rather than $\omega$.

\begin{remark}
The term `$p$-adic integration' refers to integration theory on $F$-analytic manifolds where $F$ is an arbitrary non-archimedean local field, not necessarily isomorphic to $\Qb_p$. This could be considered an abuse of language, but it's widespread and harmless.
\end{remark}

\subsection{Weil's canonical measure}\label{sub:canonical}

Let $X$ be a smooth $F$-variety (that is, a smooth separated $F$-scheme of finite type). The set of $F$-rational points $X(F)$ is endowed with a structure of an $F$-analytic manifold. This is shown as in the case of real manifolds, by using an analogue of the inverse function theorem for $F$-analytic functions. We refer the reader to \cite{MR1743467} for more details. Similarly, a regular top-degree form on $X$ induces an analytic top-degree form on $X(F)$.

In the case of a smooth $\Oo_F$-scheme $X$, the $F$-analytic manifold $X(\Oo_F)$ is endowed with a canonical density. 
\begin{theorem}[Weil]\label{thm:weil}
Let $X$ be a smooth $\Oo_F$-scheme. There exists a unique Borel measure $\mu_{X}$ on $X(\Oo_F)$, such that one has $$\mu_X|_{U(\Oo_F)} = |\omega|$$
for every $\omega \in \Omega^{\top}_{X/\Oo_F}(U)$ satisfying
$$\Omega^{\top}_{U/\Oo_F} = \Oo_U \cdot{} \omega.$$
Furthermore, the volume of $X(\Oo_F)$ can be computed in terms of the number of elements of $X(k_F)$ as follows:
\begin{equation}\label{eqn:vol}\vol_{\mu}\left(X(\Oo_F)\right) = \frac{|X(k_F)|}{q^{d}},\end{equation}
where $d$ denotes the dimension of the $k_F$-variety $X \times_{\Oo_F} k_F$.
\end{theorem}

\begin{proof}
We refer the reader to the proof of Theorem \ref{thm:formal_weil} below, respectively the references \cite{MR1743467} and \cite{weil2012adeles}.
\end{proof}

In our work we work in a more general setting than the one of smooth $\Oo_F$-varieties (see the theorem below). In the following we denote by $\Spf A$ the affine formal scheme associated to a complete local ring. We refer the reader to \cite{ALONSOTARRIO20091373} for definition and properties of smooth morphisms of formal schemes.

\begin{definition}\label{definition:formal}
Let $\Xc$ be a Noetherian formal scheme and $\phi\colon \Xc \to \Spf \Oo_F$ a smooth and separated morphism of formal schemes. We denote by $\Xc(\Oo_F)$ the set of sections $\Spf \Oo_F \to \Xc$ of $\phi$, and similarly $\Xc(k_F)$ for the set of morphisms $\Spec k_F \to \Xc$ making the following diagram commute.
\begin{equation}
\begin{tikzcd}
\ & \Xc \arrow[d, "\phi"] \\
\Spec k_F \rar \arrow[ur] & \Spf \Oo_F
\end{tikzcd}
\end{equation}
\end{definition}
Since $\Xc$ is Noetherian, $\Xc(k_F)$ is a finite set. There is a specialisation map $\sp_{\Xc}\colon \Xc(\Oo_F) \to \Xc(k_F)$, obtained by precomposing a section $\Spf \Oo_F \to \Xc$ with $\Spec k_F \to \Spf \Oo_F$. In particular we may express the set $\Xc(\Oo_F)$ as a finite disjoint union
$$\Xc(\Oo_F) = \bigsqcup_{\bar{x} \in \Xc(k_F)} \sp_{\Xc}^{-1}(\bar{x}).$$
We now define the structure of an $F$-analytic manifold on $\Xc(\Oo_F)$. To do so, we choose for every $\bar{x} \in \Xc(k_F)$ a morphism of $\Spf \Oo_F$-formal schemes
\begin{equation}\label{star}
\alpha_{\bar{x}}\colon\Spf \Oo_F[[x_1,\dots,x_{n_{\bar{x}}}]] \to \Xc, \text{which is formally \'etale and sends $0$ to $\bar{x}$},
\end{equation}
where we denote by $n_{\bar{x}}$ the rank of the sheaf of formal K\"ahler differentials $\hat\Omega^1_{\Xc/\Spf\Oo_F}$ at $\bar{x}$. 

Since $\alpha_{\bar{x}}$ is formally \'etale, it induces an isomorphism 
$$  \sp_{\Xc}^{-1}(\bar{x}) \simeq  \sp^{-1}_{\Spf \Oo_F[[x_1,\dots,x_{n_{\bar{x}}}]]}(0).$$
The right-hand side is by definition equivalent to the set of continuous $\Oo_F$-linear ring homomorphisms
$$\Oo_F[[x_1,\dots,x_{n_{\bar{x}}}]] \to \Oo_F,$$
such that the composition $\Oo_F[[x_1,\dots,x_{n_{\bar{x}}}]] \to \Oo_F \to k_F$ is the canonical map sending the variables $x_i$ to $0$.
The set of such maps is in canonical bijection with $\m_F^{n_{\bar{x}}}$ (the images of the topological generators $x_i$ give rise to this identification). The latter is an open subset of $F^{n_{\bar{x}}}$. We thus have a disjoint covering of $\Xc(\Oo_F)$ by such sets. We give $\Xc(\Oo_F)$ the induced topology. There is a unique maximal $F$-atlas on $\Xc(\Oo_F)$ containing the above charts; this makes $\Xc(\Oo_F)$ into an $F$-manifold. 
\begin{lemma}
The structure of an $F$-analytic manifold on $\Xc(\Oo_F)$ is independent of the choice \eqref{star}.
\end{lemma}

\begin{proof}
Two choices $\alpha_{\bar{x}}$ and $\beta_{\bar{x}}$ in \eqref{star} differ by a continuous $\Oo_F$-linear automorphism $\gamma$ of $$\Oo_F[[x_1,\dots,x_{n_{\bar{x}}}]].$$ Such an automorphism is represented by a formal power series 
$$\sum a_{\underline{i} \in \mathbb{N}^{n_{\bar{x}}}}(x_1,\dots,x_{n_{\bar{x}}}) \underline{x}^{\underline{i}}$$
in $x_1,\dots,x_{n_{\bar{x}}}$, where $a_{\underline{i}} \in \Oo_F[[x_1,\dots,x_{n_{\bar{x}}}]]$. Since an element of $\Oo_F$ has norm $\leq 1$ and elements of $\m_F$ have norm $< 1$, these power series converge when the variables are evaluated in $\m_F$. This implies that $\gamma$ induces an analytic automorphism of $\m_F^{n_{\bar{x}}}$, and therefore the structure of an analytic manifold is indeed independent of choices.
\end{proof}

We denote by $\widehat{\Omega}^\top_{\Xc/\Oo_F}$ the top exterior power of the sheaf of formal $\Oo_F$-linear K\"ahler differentials on $\Xc$. By definition, it is an invertible sheaf on $\Xc$.

For every isomorphism as in \eqref{star} a section $\omega$ of this sheaf gives rise to a continuous top degree form
$$f \cdot{} dx_1\wedge \cdots \wedge dx_{n_{\bar{x}}},$$
where $f \in \Oo_F[[x_1,\dots,x_{n_{\bar{x}}}]]$. As above, $f$ is a convergent power series on $\m_F^{n_{\bar{x}}}$. We therefore get an analytic top degree form $\omega$ on $\Xc(\Oo_F)$. The integral $\int_{\Xc(\Oo_F)}|\omega|$ is therefore well-defined.

We can now state the following mild generalisation of Weil's theorem.

\begin{theorem}\label{thm:formal_weil}
Let $\Xc$ be as in Definition \ref{definition:formal}. There exists a unique Borel measure $\mu_{\Xc}$ on $\Xc(\Oo_F)$, such that one has $$\mu_{\Xc}|_{\Uc(\Oo_F)} = |\omega|$$
for every $\omega \in \widehat{\Omega}^{\top}_{\Xc/\Oo_F}(\Uc)$ satisfying
$$\Omega^{\top}_{\Uc/\Oo_F} = \Oo_{\Uc} \cdot{} \omega.$$
Furthermore, the volume of $\Xc(\Oo_F)$ can be computed in terms of the number of elements of $\Xc(k_F)$ by using the following identity
\begin{equation}\label{eqn:vol}\vol_{\mu}\left(\sp^{-1}(\bar{x})\right) = \frac{1}{q^{n_{\bar{x}}}},\end{equation}
where $n_{\bar{x}}$ is the positive integer of \eqref{star}.
\end{theorem}

\begin{proof}
Let $\Uc \subset \Xc$ be an open subset, and assume that for $i=1,2$ we have sections $\omega_i \in \widehat{\Omega}^{\top}_{\Uc/\Oo_F}$, such that 
$$\widehat{\Omega}^{\top}_{\Uc/\Oo_F} = \Oo_{\Uc} \cdot{} \omega.$$
This implies $\omega_2 = f\cdot{} \omega_1$, where $f\in \Oo_{\Uc}^{\times}$ is an invertible regular function. We denote by $f\colon \Uc \to \G_{m,\Spf \Oo_F} = \Spf \Oo_F[t,t^{-1}]$ the corresponding map of formal schemes. For every $x \in \Uc(\Oo_F)$ we obtain 
$$f(x) \in \Oo_F \text{ and }f^{-1}(x) \in \Oo_F.$$
This yields $f(x) \in \Oo_F^{\times}$, and thus $|f(x)| = 1$. We conclude $|\omega_1|=|\omega_2|$. Therefore, $|\omega_1| = |\omega_2|$ and we see that the measure is well-defined and independent of the choice of the $\Oo_F$-generator $\omega$.

In particular, for every choice of an isomorphism of formal schemes \eqref{star} we have 
$$\int_{\sp^{-1}(\bar{x})} |\omega| = \int_{\m_F^{n_{\bar{x}}}} |dx_1 \wedge \cdots \wedge dx_{n_{\bar{x}}}| = \frac{1}{q^{n_{\bar{x}}}}.$$
This concludes the proof of the theorem.
\end{proof}

\subsection{$p$-adic volumes of abelian varieties}

Let $A$ be an abelian $F$-variety, that is, a smooth proper and separated group scheme defined over $\Spec F$. According to a theorem by N\'eron there exists a unique (up to a unique isomorphism) $\Oo_F$-scheme $\Nn$ over $\Spec \Oo_F$, such that the following properties hold:
\begin{enumerate}
\item[(a)] there is an isomorphism of $F$-schemes $\Nn \times_{\Oo_F} F \simeq A$,

\item[(b)] $\Nn$ is a smooth and separated $\Oo_F$-group scheme of finite type,

\item[(c)] for every smooth $\Oo_F$-scheme $S \to \Spec \Oo_F$, and an $F$-morphism
$$S \times_{\Oo_F} F \to A,$$
there is a unique extension to an $\Oo_F$-morphism
$$S \to \Nn.$$
\end{enumerate}
If $\Nn$ happens to be an abelian $\Oo_F$-scheme, $A$ is said to be of \emph{good reduction}.

Condition (c) above applies to $\Spec \Oo_F$ and amounts to an equality
\begin{equation}\label{eqn:2}
A(F) = \Nn(\Oo_F).
\end{equation}
Hence, every $F$-rational point of $A$ extends to a unique $\Oo_F$-rational point of $\Nn$.

\begin{definition}\label{defi:nor_volume}
Let $A$ be an abelian $F$-variety with N\'eron model $\Nn$. We define the \emph{normalised volume} $\tilde{V}$ of $A(F)$ to be 
$$\tilde{V}(A)=\vol_{\mu_{\Nn}}\left(\Nn(\Oo_F)\right),$$
where $\mu_{\Nn}$ refers to Weil's measure defined in Theorem \ref{thm:weil}.
\end{definition}

By virtue of Formula \eqref{eqn:vol} we have
$$\tilde{V}(A(F))=\frac{|\Nn(k_F)|}{q^d}.$$
The right-hand side of this equality can be further simplified to 
$$\frac{|\pi_0(\Nn)(k_F)|\cdot{}|\Nn^{\circ}(k_F)|}{q^d},$$
where $\pi_0(\Nn)$ denotes the \'etale sheaf on $\Spec k_F$ of relative connected components of the base change $\Nn_{k_F} \to \Spec k$, and $\Nn^{\circ} \subset \Nn$ denotes the neutral connected component. A connected commutative group scheme over a field, is isomorphic to an extension of an abelian variety by an affine commutative group scheme (\emph{Chevalley d\'evissage}). The latter is an extension of a torus by a unipotent group scheme. 

Under the assumption that the toric part is split (which will always be satisfied in our work), one obtains 
$$\tilde{V}(A(F))=\frac{|\pi_0(\Nn)(k_F)|\cdot{}q^{r_u}\cdot{}(1-q)^{r_t}}{q^d},$$
where $r_u$ and $r_t$ denote the dimensions of the unipotent, respectively toric, part of $\Nn^{\circ} \times_{\Oo_F} k_F$.

\begin{lemma}
Let $A$ be an abelian $F$-variety, and $\omega \in \Omega^{\top}_{A/F}$ a regular top-degree form on $A$. We introduce the notation $V_{\omega}(A)$ to denote $\int_{A(F)} |\omega|$. There exists a positive rational $c_{\omega}(A)$, called the \emph{conductor}, such that
$$V_{\omega}(A) = c_{\omega}(A)\cdot{} \tilde{V}(A).$$
\end{lemma}

\begin{proof}
Let $\Nn \to \Spec \Oo_F$ be a N\'eron model for $A$. We view $\omega$ as a rational section of the invertible sheaf $\Omega^{\top}_{\Nn/\Oo_F}$.

Since $\omega$ is a regular top-degree form on $A$, it is translation-invariant. This implies that the pole order $\nu$ of $\omega$ along the divisor $\Nn \times_{\Oo_F} k_F \subset \Nn$ is constant along this (possibly disconnected) divisor. Let $\unif \in F$ be an element of valuation $1$. We then have that
$$\unif^{\nu} \cdot{} \omega$$
is a regular section of $\Omega^{\top}_{\Nn/\Oo_F}$, and furthermore satisfies
$$\Omega^{\top}_{\Nn/\Oo_F} = \Oo_{\Nn} \cdot{} (\unif^{\nu}\cdot{}\omega).$$
We conclude from Weil's Theorem \ref{thm:weil} that $\tilde{V}$ is given by integrating the density $|\unif^{\nu}\cdot \omega| = |q|^{\nu}\cdot |\omega|$. Therefore one has 
$V_{\omega}(A) = q^{-\nu}\cdot{}\tilde{V}(A)$.
\end{proof}

\subsection{A determinantal formula }

Recall that $A$ is the ring $\Z[[T_e : e \in E]]$. Let $B(\Gamma)$ be the affine scheme $\Spec \cO_F[[T_e : e \in E]]$; it carries a natural map to $\Spec A$. Let $\Delta$ be the subscheme given by the union of divisors $(T_e)$. Consider the intersection $B(\Gamma)^{\flat} := B(\Gamma)(\cO_F) \cap (B(\Gamma) \setminus \Delta)(F)$. 

We have a scheme $\pversal(\Gamma, \eta)$ over $\Spec A$ (see Definition \ref{def:versalspace}) for each $\eta \in \redC_0(\Gamma, \Zb)$. Given $a \in B(\Gamma)^{\flat}$, the fiber 
$$\pversal(\Gamma, \eta) \times_{\Spec A} a$$
is an abelian variety. We will study its normalised volume. 
\begin{lemma} \label{periodmap}
Let $\periodmap_a : H_1(\Gamma, \Z) \to H^1(\Gamma, F^{\times})$ be the homomorphism defined by composing $\sigma$ with the evaluation $H^1(\Gamma, \G_{m, \Kk_{\Gamma}}) \to H^1(\Gamma, F^{\times})$ induced by $\pi_e \to a_e$. We have
\[ \pversal(\Gamma, \eta) \times_{\Spec A} a =  H^1(\Gamma, F^{\times}) / \periodmap_a(H_1(\Gamma, \Z)) \]
\end{lemma}
\begin{proof}
This is a direct application of the results of \cite[\S 3.8]{MR1083353} to the period lattice described by Lemma \ref{lem:genericfiber}.
\end{proof}

The volume will be expressed in terms of the following `logarithm' of $\sigma_a$. Let $\{ x_e \}$ for $e \in E(\Gamma)$ be a tuple of integers, and consider the mapping $\ctc_x \colon H_1(\Gamma, \Z) \to H^1(\Gamma, \Z)$ defined by
\begin{equation} \label{definition:ctc} \gamma \to \sum_{e \in E(\Gamma)} x_e \langle \gamma, e \rangle [e].\end{equation}
Note that this map does not depend on the orientation of $\Gamma$.

\begin{proposition} \label{proposition:ctcinjective}
	Suppose $x_e > 0$ for all edges. Then the map $\ctc_x$ is injective. 
\end{proposition}

\begin{proof}
	Suppose $\gamma \neq 0$. We have $\langle \gamma, \ctc_x(\gamma) \rangle = \sum_{e \in E(\Gamma)} x_e \langle \gamma, e \rangle^2 > 0$.
\end{proof}

For $a \in \pbase(\Gamma)^{\flat}$, we have a tuple of integers $\{\nu_e(a) \}$ which we abbreviate to $\nu(a)$. 
\begin{proposition} \label{prop:columeoffibre}
\begin{equation}\label{equation:volumeoffibre} \tilde{V}(\fibera) =  \left( \frac{q-1}{q} \right)^{h^1(\Gamma)} |H^1(\Gamma,\Z) / \ctc_{\nu(a)}(H_1(\Gamma,\Z))|. \end{equation}
\end{proposition}
\begin{proof} 
By virtue of Definition \ref{defi:nor_volume} we have 
$$\tilde{V}(\fibera) = q^{\dim \fibera} \cdot{} |\Nn^{\circ}(k_F)| \cdot{} |\Phi(k_F)|,$$
where $\Nn/\Oo_F$ denotes the N\'eron model of $\fibera/F$, $\Nn^{\circ}$ the neutral connected component and $\Phi$ the \'etale $k_F$-group scheme of connected components of the special fibre of the N\'eron model.

By construction, $\fibera$ has split multiplicative reduction, thus $\Nn^{\circ}_{k_F}$ is isomorphic to a split torus $\G_{m,k_F}^{h^1(\Gamma)}$. This explains the factor $(1-q^{-1})^{h^1(\Gamma)}$ It remains to verify that the Tamagawa number $|\Phi|$ is equal to $|H^1(\Gamma,\Z) / \ctc_{\nu(a)}(H_1(\Gamma,\Z))|$. This follows from a description of Tamagawa numbers of abelian varieties $A/F^{un}$ with a non-archimedean uniformisation
$$M \hookrightarrow \G_{m,F^{un}}^r \twoheadrightarrow A,$$
here $M$ denotes a free abelian group, and $F^{un}/F$ denotes the maximal unramified extension of $F$. 
As explained in \cite[p. 82]{halle-nicaise}, there is an exact sequence
$$0 \to (\Zb^r/M^I)_{tors} \to \Phi(A) \to H^1(I,M),$$
where $I$ denotes the absolute Galois group of $F^{un}$, that is, $I=\mathsf{Gal}(F^s/F^{un})$.

The universal property of the N\'eron model implies that it commutes with base change along $F^{un}/F$, and therefore we are free to replace $\fibera$ by $A = \fibera \times_{\Spec F} \Spec F^{un}$. In this case we have $M=H_1(\Gamma,\Zb) = \Zb^r$ with the trivial Galois action, and therefore
$H^1(I,\Zb^r) = 0$ (this follows from the long exact sequence associated to $\Zb \hookrightarrow \Qb \twoheadrightarrow \Qb/\Zb$).

Applying this to the non-archimedean uniformisation
$$H_1(\Gamma,\Zb) \hookrightarrow H^1(\Gamma,\G_{m,F^{un}}) \twoheadrightarrow A,$$
we obtain $\Phi(A) \simeq H^1(\Gamma,\Zb)/\tau_{\nu(a)}(H_1(\Gamma,\Zb))$ as asserted above.
\end{proof}

\subsection{The volume polynomial is the Kirchhoff polynomial}
We may view $|H_1(\Gamma, \Z) / \ctc_{x}(H_1(\Gamma, \Z))|$ as a function of the variables $x_e \geq 0, e \in \Gamma$. Pick dual bases of $H^1(\Gamma, \Z)$ and $H_1(\Gamma,\Z)$. The mapping $\ctc_\Gamma$ defines an $r \times r$ matrix $T_{\Gamma}$ with entries $\langle \gamma_i, \ctc_x(\gamma_j) \rangle$. The proof of Proposition \ref{proposition:ctcinjective} shows that it is positive definite when $x_e > 0$. By standard linear algebra, we obtain the following.
\begin{lemma}
For $x_e > 0$ we have $|H_1(\Gamma, \Z) / \ctc_{x}(H_1(\Gamma, \Z))| = \det T_{\Gamma}$.
\end{lemma} 
In particular, the right-hand side is a homogenous element of $\Zb[(x_e)_{e \in E(\Gamma)}]$ of total degree $h_1(\Gamma)$. 

A subset $T \subset \edgesgam$ is said to be a \emph{spanning tree} if it is connected an simply connected, and is incident to all vertices of $\Gamma$. A graph has a spanning tree if and only if it is connected. In the disconnected case we choose a spanning tree for every connected component and refer to the resulting subgraph as a \emph{maximal spanning forest}.

\begin{definition}
Let $\Arb(\Gamma)$ be the set of maximal spanning forests of $\Gamma$. Define the Kirchhoff polynomial of $\Gamma$ by
$$\Psi_{\Gamma}=\sum_{T \in \Arb(\Gamma)} \prod_{e \notin T} x_e.$$
\end{definition}
The function $\Psi_{\Gamma}$ is therefore the basis generating function for the bond matroid of $\Gamma$. Basis generating functions of matroids satisfy several interesting properties. See Br\"and\'en--Huh \cite{bh} for a modern viewpoint through the lens of Lorentzian polynomials.

\begin{proposition}
\[ \det T_\Gamma = \Psi_{\Gamma} \]
\end{proposition}
\begin{proof}
The determinant $\det T_\Gamma$ is a sum of monomials of degree $h_1(\Gamma)$. It is thus enough to understand its specialisations at $x_e = 0$ for $e \notin S$, where $S$ is a subset of $\Gamma$ with $|S| \leq h_1(\Gamma)$. If the rank of $T_{\Gamma}$ under this specialisation is less than $h_1(\Gamma)$, then the determinant vanishes. Thus we may assume that $|S| = h_1(\Gamma)$ and  that $[e], e \in S$ are linearly independent in $H^1(\Gamma, \Q)$. Any such set of edges forms a basis in $H^1(\Gamma, \Z)$. Picking a dual basis for $H_1(\Gamma, \Z)$, we see that $ \det T_{\Gamma} = \prod_{e \in S} x_e$ in this specialisation. Since $\Gamma \setminus S$ is a spanning forest exactly when $[e], e \in S$ form a basis of $H_1(\Gamma, \Z)$, the result follows. \end{proof}

\begin{corollary}[Main result]\label{cor:main}
For $a \in \base(\Gamma)(\Oo_F) \setminus \Delta(\Oo_F)$ we have
$$\vol (\hmap^{-1}(a)) = (1-q^{-1})^{h_1(\Gamma)}\cdot{} \Psi_{\Gamma}(\nu(a)).$$
\end{corollary}

Recall that an edge $e$ is said to be a \emph{loop} if it is incident to only one vertex. We call $e$ a bridge, if $\Gamma \setminus e$ is disconnected. An edge which is neither a loop nor a bridge is said to be \emph{ordinary}. We omit the elementary proof of the following lemma.

\begin{lemma}[Deletion-contraction]\label{lemma:deletion-contraction}
Let $e$ be an ordinary edge, then
\begin{equation}\label{eqn:ordinary}
\Psi_{\Gamma} = x_e\cdot{}\Psi_{\Gamma \setminus e} + \Psi_{\Gamma/e}.
\end{equation}
If $e$ is a loop, then 
\begin{equation}\label{eqn:loop}
\Psi_{\Gamma} = x_e \cdot{}\Psi_{\Gamma \setminus e}.
\end{equation} 
If $e$ is a bridge, then 
\begin{equation}\label{eqn:bridge}
\Psi_{\Gamma} = \Psi_{\Gamma/e}.
\end{equation}
Furthermore, the map $\Gamma \mapsto \Psi_{\Gamma}$ is the unique map satisfying the relations above, and $\Psi_{\Gamma_0} = 1$ for a graph $\Gamma_0$ with an empty edge-set $E(\Gamma_0)$.
\end{lemma}

We now give a second, closely related interpretation of $\Psi_{\Gamma}$. Let $\underline{n}=(n_e)_{e\in E(\Gamma)}$ be an tuple of positive integers indexed by $E(\Gamma)$. We denote by $\Gamma_{\underline{n}}$ the graph obtained by replacing $e$ by a linear tree with $n$ vertices. We refer to this linear tree as the $n_e$-fragmentation of $e$, and to $\Gamma_{\underline{n}}$ as the \emph{\underline{n}-fragmentation} of $\Gamma$. The lemma above implies the following:

\begin{corollary}
The number of maximal spanning forests in the fragmentation $\Gamma_{\underline{n}}$ is equal to the specialisation $\Psi_{\Gamma}|_{x_e=n_e}$.
\end{corollary}

In light of this corollary we obtain another interpretation of our main result Corollary \ref{cor:main}. We see that the $p$-adic volume is related to the number of spanning trees in fragmentations of $\Gamma$. These fragmented graphs have also appeared in computations of motivic zeta functions of Jacobians by Halle--Nicaise, see the proof of Lemma 4.3.1.3 in \cite{halle-nicaise}.

\subsection{Integrating along the integrable system}

In this section we explain how $p$-adic integration can be used to count the number of $k_F$-rational points of $\hmap^{-1}(0)$. This method was used by one of the authors in joint work with Wyss and Ziegler in \cite{gwz}. In particular, the definition of the top degree form $\omega$ on $\pDol$ mimics the constructions of \emph{loc. cit.}

We emphasise that in the present situation it would not be difficult to directly count the number of $k_F$-rational points with different methods. However, we believe that the calculations below are sufficiently pleasant to be included as a short reality check.

The space $\pDol=\pDol(\Gamma)$ contains an open subset $\pDol' \subset \pDol(\Gamma)$ that is a (trivial torsor under a) commutative group scheme $\mathcal{P} \to \base$. For $a \in \base \setminus \Delta$ one has 
$$\pDol' \times_{\base} a = \pDol \times_{\base} a,$$
and furthermore, $\pDol'$ is surjective over $\base$. This shows that $\pDol'$ has a complement of codimension $\geq 2$.

Therefore, the sheaf of relative top degree forms $\Omega^{\top}_{\pDol'/\base}$ can be trivialised. Let $\omega_0 \in \Omega^{\top}_{\pDol'/\base}(\pDol')$ be a generating section. We denote by $\omega_1 \in \Omega^{\top}(\base)$ a generating section of formal top degree forms on $\base$. The wedge product $\omega = \omega_0 \wedge \omega_1$ is therefore a rational top degree form on $\pDol$ defined up to codimension $2$. By virtue of Hartogs's Extension Theorem we have that 
$\omega$ is in fact everywhere well-defined, and a generating section of $\Omega^{\top}_{\pDol}$.

Using Fubini we can factor the integral of $|\omega|$ along the map  $\hmap : \pDol(\Gamma) \to \pcyclebase(\Gamma)$ to obtain 
\begin{equation} \label{equation:integrateoverthebase} V_{\omega}(\pDol(\Gamma)) :=  \int_{a \in \m_F^{|E(\Gamma)|}\cap H_1(\Gamma)} V_{\omega_0}(\fibera(F)) |\omega_1|.  \end{equation}

\begin{lemma}
For $a \in \base(\Gamma)^{\flat}$ one has has $V_{\omega_0}(\fibera(F)) = \tilde{V}(\fibera(F))$.
\end{lemma}

\begin{proof}
We consider the base change
$$\cM= \pversal(\Gamma,\eta)' \times_{\base,a} \Spec \Oo_F.$$
The universal property of the N\'eron model of $\fibera$ yields a map $\cM \to \cN$. Since $\cM$ is a group scheme on $\Oo_F$ that has semi-abelian reduction, this map is an open immersion. The restriction of the top degree form $\omega_0$ to $\cM$ is still a generating section. It extends uniquely to a translation-invariant section of $\cN$ which is still a generating section (by translation-invariance). This concludes the proof of the identity $V_{\omega}(\fibera(F)) = \tilde{V}(\fibera(F))$.
\end{proof}

Let $\nu \in \N^{E(\Gamma)}$ and let $\pcyclebase(\Gamma)^{\nu} \subset \pcyclebase(\Gamma)$ be the subset of points with valuation $\nu$. The fibres of $\hmap : \pDol(\Gamma) \to \pcyclebase(\Gamma)$ over a fixed $\pcyclebase(\Gamma)^{\nu}$ have constant volume given by Proposition \ref{prop:columeoffibre}. Thus the right-hand side of \eqref{equation:integrateoverthebase} equals
\begin{equation}  \label{equation:totalvolumesum} \left( \frac{q-1}{q} \right)^{h^1(\Gamma)} \sum_{\nu \in \N^{E(\Gamma)}} \Psi_{\Gamma}(\nu) \mu(\pcyclebase(\Gamma)^\nu). \end{equation}
We now illustrate this decomposition in a few simple examples. Let $F(\Gamma) = q^{h^1(\Gamma)} V_\omega(\pDol(\Gamma))$. 

\begin{example}
	Let $\boing$ be the graph with one vertex and one edge. Then 
	\[ F(\boing) = (q-1) \sum_{n \in \N^{>0}} n (q-1) q^{-n-1} = \frac{(q-1)^2}{q} \sum_{n \in \N^{>0}} n q^{-n} =  \frac{(q-1)^2}{q}   \frac{q^{-1}}{(1-q^{-1})^2} = 1 \]
\end{example} 
\begin{example}
	More generally, if $\Gamma$ is a loop with $N$ vertices and $N$ edges, we have
	\[  F(\Gamma) = (q-1) \sum_{n \in \N^{>0}} Nn q^{-n-1} = N. \]
\end{example}

In more general cases, the decomposition along the base $\pcyclebase$ becomes more complicated. However, we can prove the following recursion for the volume by working "uniformly" on the fibres.

\begin{proposition} \label{proposition:totalvolumedeletioncontraction} For an ordinary edge $e$ we have
	\begin{equation}\label{eqn:F} F(\Gamma) = F(\Gamma / e) + F(\Gamma \setminus e). \end{equation}
	For a bridge $e$ we get $F(\Gamma) = F(\Gamma / e)$, and for a loop $e$, we get $F(\Gamma) = F(\Gamma \setminus e)$.
\end{proposition}
\begin{proof}
	Assume by induction that the formula has already been established for graphs with a smaller number of edges than $\Gamma$. Suppose that $\Gamma$ contains an ordinary edge $e$. 
	As we have seen, the function $\tilde{V}(\Gamma)$ equals $\frac{q-1}{q}\cdot{}\nu_e \cdot{} \tilde{V}(\Gamma \setminus e) + \tilde{V}(\Gamma/e)$ for every ordinary edge $e$. Integrating along the $e$-variable over $\m_F$, we obtain
	$$\int_{\m_F} \tilde{V}(\Gamma)|dx_e| = \int_{\m_F} \nu_e \cdot{} \tilde{V}(\Gamma \setminus e) |dx_e| + \int_{\m_F} \tilde{V}(\Gamma / e) |dx_e|.$$
	This yields
	$$\int_{\m_F} \tilde{V}(\Gamma)|dx_e| = \frac{1}{q}\cdot{}\left(\frac{q}{q-1}\right)\cdot{} \left(\frac{q-1}{q}\right) \cdot{} \tilde{V}(\Gamma \setminus e) + \int_{\m_F} \tilde{V}(\Gamma / e) |dx_e|.$$
	Note that the rightmost integral is not further simplified, since the function $\tilde{V}(\Gamma/e)$ may still depend on $x_e$. This can happen if there is an edge $e'$, satisfying $e \neq e'$ and $[e] = [e']$.
	
	Integrating over the remaining variables and multiplying with $q^{h^1(\Gamma)}$ we obtain
	$$F(\Gamma) = q^{h^1(\Gamma)}\cdot{} q^{-1}\cdot \int_{\m_F^{|E(\Gamma) \setminus e|} \cap H_1(\Gamma \setminus e)}\tilde{V}(\Gamma)|\wedge_{e' \neq e}dx_{e'}| + q^{h^1(\Gamma)}  \int_{\m_F^{|E(\Gamma) \setminus e|}\cap H_1(\Gamma / e)} \tilde{V}(\Gamma/e).$$
	Since $h^1(\Gamma) = h^1(\Gamma/e)$, the right-hand side agrees with the sum $F(\Gamma \setminus e) + F(\Gamma / e)$.
	This concludes the proof if $\Gamma$ contains an ordinary edge. The remaining cases are proven by similar computations which we leave to the reader.
\end{proof}
\begin{corollary} \label{proposition:volumecountsforests}
	$F(\Gamma) = \Psi_{\Gamma}(\vec{1}) = \#\{\text{maximal forests in $\Gamma$}\}$.
\end{corollary}
\begin{proof}
	By Proposition \ref{proposition:totalvolumedeletioncontraction} and Lemma \ref{lemma:deletion-contraction}, both $F(\Gamma)$ and $\Psi_\Gamma(\vec{1})$ satisfy the same recursion. Both quantities are equal to $1$ for the base case of this recursion, i.e. a graph with no edges.
\end{proof}

\begin{corollary} \label{proposition:centralfibresise}
	\[	|\pDol(\Gamma)(k_F)| = \#\{\text{maximal forests in }\Gamma \} \times q^{h_1(\Gamma)}.  \]
\end{corollary}
\begin{proof}
	Combine Theorem \ref{thm:formal_weil} and Corollary \ref{proposition:volumecountsforests}.
\end{proof}

\section{Tropicalisation}\label{tropical}

Henceforth we will use the term \emph{local field} to designate non-archimedean and archimedean examples.

\subsection{Tropical abelian varieties and normalised volume}

Let $K$ be a normed field. The \emph{tropicalisation} of the algebraic $K$-torus $\G_{m,K}^r$ is defined to be the set
$$\trop (\G_{m,K}) = \Rb^r.$$
There is a map (of sets) $(K^{\times})^r \to \Rb^r$ given by $\log |\cdot{}|$. For a subvariety $Z \subset \G_{m,K}^r$, the tropicalisation is defined as 
$$\trop(Z) = \overline{\log |Z(\bar{K})|} \subset \Rb^r.$$

Let $K$ be a local field, and $A/K$ an abelian variety which can expressed as a quotient $\G_{m,K}^r/Y$, where $Y$ is a discrete group. Then, one defines
$$\trop(A) = \trop(\G_{m,K}^r)/\trop(Y) = \Rb^r/\log|Y|.$$
The latter is an $r$-dimensional topological torus. 

The Lebesgue measure on $\Rb^r$ descends to a measure on the quotient $\trop(A)$ that we will denote by $\mu_{\Rb}$. 

\begin{proposition}\label{prop:trop}
Let $K$ be a non-archimedean local field. Then, we have $\tilde{V}(A(K)) = (1-q^{-1})^{\dim_K A}\cdot{} \vol_{\mu_{\Rb}}(\trop(A))$.
\end{proposition}

\begin{proof}
We choose the same notation as before, and write $A=\G_{m,K}^r/Y$. Since $A(K) = (K^{\times})^r/Y$, and $K^{\times} = \Zb \times \Oo_K$, one has
$$\tilde{V}(A) = |\Zb^r / \nu(Y)|\cdot{} q^{-r}.$$ The number of points in the quotient is equal to the absolute value of the determinant of the map $\Zb^r \to \Zb^r$ whose image agrees with $Y$. This concludes the proof, since the same expression computes the volume of the cokernel $\mathsf{coker}(\Zb^r \to \Rb^r) = \trop(A)$.
\end{proof}

\subsection{Hypertoric Hitchin fibres and tropical Jacobians}

A \emph{tropical curve} is a pair $(\Gamma,w)$ where $\Gamma$ is a graph and $w\colon E(\Gamma) \to \Qb_{> 0}$ is a function. Its Jacobian is defined to be the cokernel
$$J(\Gamma,w) = \mathsf{coker}(\tau_w),$$
where we recall that $\tau_w$ denotes the map $H_1(\Gamma,\Zb) \to H^1(\Gamma,\Zb)$ given by 
$$\gamma \mapsto \sum_{e \in E(\Gamma)} w(e) \langle \gamma, [e]\rangle [e].$$

The relationship between the tropical Jacobian of the tropicalisation of a smooth curve $C/F$ and the tropicalisation of $J(C)/F$ is as clear-cut as it could be: by virtue of a theorem of Baker--Rabinoff \cite{baker-rabinoff} the two are isomorphic. 

\begin{proposition}
Let $F$ be a non-archimedean local field. For every $b \in B(\Gamma)^{\flat}$ there is a tropical curve $(\Gamma_b,w_b)$, such that there is a measure-preserving isomorphism
$\trop(\pi^{-1}(b)) \simeq J(\Gamma_b,w_b).$
\end{proposition}

\begin{proof}
Let $\Gamma_b = \Gamma$ and $w_b(e) = \nu_e(b)$ for all $e \in E(\Gamma)$. According to \eqref{periodmap} ,
$$\pi^{-1}(b) = H^1(\Gamma,\G_{m,K}) / \sigma_b(H_1(\Gamma,\Zb)).$$
Tropicalisation therefore yields
$$\trop(\pi^{-1}(b)) = H^1(\Gamma,\Rb)/\tau_b(H_1(\Gamma,\Zb))=J(\Gamma,w_b).$$
This concludes the proof.
\end{proof}

\subsection{Tropicalising the total space}

By virtue of the previous subsection, the set $B(\Gamma)^{\flat}$ parametrises a family of tropical curves, such that the tropicalisation of the fibre over $b \in B(\Gamma)^{\flat}$ agrees with the corresponding tropical Jacobian. We conclude our excursion into tropical geometry by extending the tropicalisation procedure described above to the total spaces. That is, we will define a morphism of manifolds with corners endowed with a measure
$$\trop(\pi)\colon \trop(\pversal(\Gamma,\eta)) \to \trop(B(\Gamma)),$$
such that the generic fibres (with respect to the fibrewise measure) are isomorphic to the tropicalisation of the corresponding hypertoric Hitchin fibre.

\begin{definition}
Let $(M,\omega)$ be a symplectic manifold which is toric with respect to the action of a torus $(S^1)^r$. We denote the quotient $M/(S^1)^r$ by $\mathsf{Del}(M)$ and refer to it as the \emph{Delzant space} of $M$.
\end{definition}

According to Karshon--Lerman's \cite[Proposition 1.1]{karshon-lerman}, the Delzant space is a manifold with corners, and the moment map gives rise to a unimodular local embedding $\mu\colon \mathsf{Del}(M) \hookrightarrow \Rb^r$. 

Recall the definition of the toric variety $\mathbb{W}_{K}$ from \ref{defi:W}. It is toric with respect to the algebraic torus $\G_{m,K}^2$. In particular, we obtain a toric symplectic variety by considering the associated complex manifold. By virtue of the definition, $\tilde{\pversal}_K(\Gamma,\eta)$ is a toric $K$-variety with dense torus isomorphic to 
$$C^1(\Gamma,\G_{m,K}^2)/C^0(\Gamma,\G_{m,K}),$$
where we embed $\G_{m,K}$ antidiagonally into $\G_{m,K}^2$. 
In the definition below we consider the associated toric symplectic manifold, obtained by setting $K=\Cb$ and passing to the world of manifolds.
\begin{definition}
\begin{enumerate}
\item[(a)] We define $\trop(\tilde{\pversal}(\Gamma,\eta)(\Cb))$ to be the quotient $\mathsf{Del}\left( \tilde{\pversal}(\Gamma,\eta)(\Cb) \right) / H_1(\Gamma,\Zb)$.
\item[(b)] The tropicalisation of the base is defined as $\trop(B(\Gamma)) = \mathsf{Del}(B(\Gamma)) = [0,1)^{E(\Gamma)}$
\item[(c)] The map $\trop(\pi)$ is defined to be the unique continuous map obtained by passing to the quotients:
$$\tilde{\pversal}(\Gamma,\eta)(\Cb)/(S^1)^{h^1(\Gamma,\Zb) + |E(\Gamma)|} \to \mathbb{C}^{E(\Gamma)}/(S^1)^{|E(\Gamma)|}.$$
\item[(d)] The tropicalisation of $\pversal(\Gamma,\eta)$ is defined to be the quotient
$$\trop(\pi)^{-1}([0,1)^{E(\Gamma)}) / H_1(\Gamma,\Zb).$$
\end{enumerate}
\end{definition}

The tropicalisations of $\pversal(\Gamma,\eta)$ and $B(\Gamma)$ have a natural piecewise linear structure. Indeed, this is obvious for the latter, since it agrees with the half-open cube $[0,1)^{|E(\Gamma)|}$. Similarly, the Delzant space of 
$\tilde{\pversal}(\Gamma,\eta)$ is isomorphic to its moment polytope, which is the source of the piecewise linear structure of $\trop(\pversal(\Gamma,\eta))$.

It would be interesting to understand the connection between the tropical Hitchin map introduced above and the tropicalisation of the universal compactified Jacobian studied in \cite{abreu-pacini,4authors}.

\bibliographystyle{amsalpha}
\bibliography{master}

\end{document}